\documentclass{jnmp}

\usepackage{amsmath}
\JNMPnumberwithin{equation}{section}

\newtheorem{theorem}{Theorem}
\newtheorem{lemma}{Lemma}
\newtheorem{proposition}{Proposition}
\newtheorem{remark}{Remark}

\theoremstyle{definition}

\begin{document}
%  Headings
%
\renewcommand{\evenhead}{D. Lukkassen and P. Wall}
\renewcommand{\oddhead}{On Weak Convergence of Locally Periodic Functions}

%  Titlepage
%
\thispagestyle{empty}

\FirstPageHead{9}{1}{2002}{\pageref{firstpage}--\pageref{lastpage}}{Article}
%  Parameters: Volume, number, year, page range, paper type
%  'Article' could be changed to 'Letter' or 'Review Article'

\copyrightnote{2002}{D Lukkassen and P Wall}
\resetfootnoterule
\Name{On Weak Convergence\\ of Locally Periodic Functions}

\label{firstpage}

\Author{Dag LUKKASSEN~$^\dag$ and Peter WALL ~$^\ddag$}

\Address{$^\dag$ Narvik Uninversity College, N-8505 Narvik, Norway \\
~~E-mail: dl@hin.no\\[10pt]
$^\ddag$ Department of Mathematics, Lule\aa\ University of Technology,
S-971 87 Lule\aa, Sweden \\
~~E-mail: walll@sm.luth.se}

\Date{Received March 03, 2001; Revised August 03, 2001; 
Accepted August 17, 2001}

\begin{abstract}
\noindent
We prove a generalization of the fact that periodic functions converge
weakly to the mean value as the oscillation increases. Some convergence
questions connected to locally periodic nonlinear boundary value 
problems are also considered.
\end{abstract}

\section{\ Introduction}

In the the proof of the reiterated homogenization results obtained in \cite
{llpw1} (see also \cite{llpw2}) the following two facts were used (the exact
defintions and properties are given later in Section 2 and 3):

\begin{itemize}
\item  If $(v_h)$ is a sequence of $Y$-periodic functions in $L_{\text{loc}%
}^p(R^n),$ $p>1,$ such that $v_h\rightarrow v$ weakly in $L^p(Y)$ and if $w_h
$ is defined as $w_h(x)=v_h(hx)$, then $w_h\rightarrow (1/\left| Y\right|
)\int_Yv(x)\,dx.$

\item  If 
\[
-\hbox{div}(a(x,hx,\xi +Du_h^\xi ))=0\text{ }
\]
on $Y=[0,y]^n,\ u_h^\xi \in W_{\text{per}}^{1,p}(Y),$ then $u_h^\xi
\rightarrow u_0^\xi $ weakly in $W_{\text{per}}^{1,p}(Y),$ where $u_0^\xi $
is the solution of the corresponding limit-equation. Here $a$ is monotone,
continuous and satisfies suitable coerciveness and growth conditions in the
third variable and periodic in the second.
\end{itemize}

We have not found proofs of these facts in the literature. The aim of this
paper is to present such proofs. Moreover, we show that the first statement
also holds for the case $p=1$.

The two facts described above are used in the proof of the reiterated
homogenization result for monotone operators, see \cite{llpw1} and \cite
{llpw2}. The solution $u_h^\xi $ is used to define a sequence of functions
similar to the ones in Tartar's celebrated method of oscillated test
functions (see e.g. the book \cite{cioranescu99}). The first fact described
above in combination with compensated compactness is used to analyze the
asymptotic behavior of this sequence of functions.

For information concerning reiterated homogenization we recommend the papers 
\cite{llpw1} and \cite{llpw2} and the references given there. Concerning
explicit engineering applications see e.g. \cite{bylund2001}.

\section{A weak convergence result}

Let us first recall the following lemma (for the proof see e.g. \cite
{dellacherie1975}).

\begin{lemma}
\label{lemdunford}\label{lemcompensated}Let $\left\{ u_h\right\} $ be a
sequence in $L^1(\Omega ).$ The following statements are equivalent:

\begin{enumerate}
\item  every subsequence of $\left\{ u_h\right\} $ contains a subsequence
which converges weakly in $L^1(\Omega ).$

\item  for all $\varepsilon >0$ there exists $t_\varepsilon >0$ such that
for all $h\in \left\{ 1,2,...\right\} $ it holds that 
\[
\int_{\left\{ u_h\geq t_\varepsilon \right\} }\left| u_h\right| dx\leq
\varepsilon ,
\]
where $\left\{ u_h\geq t_\varepsilon \right\} $ denotes the set $\left\{
x\in \Omega :u_h(x)\geq t_\varepsilon \right\} $.
\end{enumerate}
\end{lemma}

The following Proposition is a generalization of the well-known fact that a
periodic function converges weakly to its mean value as the oscillation
increases.

\begin{theorem}
\label{lemfundamental}Let $1\leq p\leq \infty $ and let $u_h\in L_{\text{loc}%
}^p(R^n)$ be $Y$ periodic for $h\in N$. Moreover, suppose that $%
u_h\rightarrow u~$weakly in $L^p(Y)$ (weak-$*$ if $p=\infty $) as $%
h\rightarrow \infty $. Let $w_h$ be defined by $w_h(x)=u_h(hx\dot{)}$. Then
as $h\rightarrow \infty $ it holds that 
\[
w_h\rightarrow \frac 1{\left| Y\right| }\int_Yu(x)dx
\]
weakly in $L^p(\Omega )$ (weak-$*$ if $p=\infty $).
\end{theorem}

%TCIMACRO{\TeXButton{Proof}{\proof}}
%BeginExpansion
\proof%
%EndExpansion
We first consider the case $1<p\leq \infty $. For simplicity we put $%
Y=(0,1)^{n}$, i.e. the unit cube in $R^{n}$, since the general case is
principally the same. Let $Y_{h}^{k}=(1/h)\left( k+Y\right) $, where $k\in
Z^{n}$, i.e. the translated image of $1/hY$ by the vector $k/h$. We note
that $\left| Y_{h}^{k}\right| =\left( 1/h\right) ^{n}$. Let $\phi \in
C_{0}(\Omega )$ and $\phi ^{h}$ the function which takes a constant value
equal to the value $\phi (k/h)$ in each cell $Y_{h}^{k}$. Due to the uniform
continuity of $\phi $ on the compact set $\overline{\Omega }$, we obtain
that $\phi ^{h}\rightarrow \phi $ uniformly on $\Omega $. Thus, 
\begin{equation}
\int_{\Omega }\left| \phi ^{h}-\phi \right| ^{q}dx\leq \int_{\Omega }\left(
\max_{x\in \Omega }\left| \phi ^{h}(x)-\phi (x)\right| \right)
^{q}dx\rightarrow 0,  \label{lqconverg}
\end{equation}
as $h\rightarrow \infty $, i.e. $\phi ^{h}\rightarrow \phi $ in $%
L^{q}(\Omega )$. Since $\phi $ has compact support in $\Omega $ we have that
each cell $Y_{h}^{k}$, for which $\phi (k/h)\neq 0$, is contained $\Omega $
for sufficiently large values of $h$. This and the $Y$-periodicity of $u_{h}$
implies that

\begin{eqnarray}
\int_{\Omega }u_{h}\left( hx\right) \phi ^{h}(x)dx &=&\int_{\mathbf{R}%
^{n}}u_{h}\left( hx\right) \phi ^{h}(x)dx  \nonumber \\
&=&\sum_{k\in \mathbf{Z}^{n}}\int_{Y_{h}^{k}}u_{h}\left( hx\right) \phi
(k/h)dx  \nonumber \\
&=&\sum_{k\in \mathbf{Z}^{n}}\left( \frac{1}{\left| Y_{h}^{k}\right| }%
\int_{Y_{h}^{k}}u_{h}\left( hx\right) dx\right) \left( \int_{Y_{h}^{k}}\phi
(k/h)dx\right)  \nonumber \\
&=&\frac{1}{\left| Y\right| }\int_{Y}u_{h}\left( x\right) dx\sum_{k\in 
\mathbf{Z}^{n}}\left( \int_{Y_{h}^{k}}\phi (k/h)dx\right)  \label{large} \\
&=&\frac{1}{\left| Y\right| }\int_{Y}u_{h}\left( x\right) dx\int_{\Omega
}\phi ^{h}(x)dx,  \nonumber
\end{eqnarray}
for sufficiently large values of $h$. Moreover, we have that $u_{h}\left(
h\cdot \right) $ is bounded in $L^{p}(\Omega )$. This fact is shown as
follows: Define the index set $I_{h}$ as 
\[
I_{h}=\left\{ k\in Z^{n}:Y_{h}^{k}\cap \Omega \neq \emptyset \right\} . 
\]
Since $\Omega $ is bounded there exists a constant $K$ such that the number
of elements in $I_{h}$ is less than $Kh^{n}$. We obtain that 
\begin{eqnarray}
\int_{\Omega }\left| u_{h}\left( hx\right) \right| ^{p}dx &\leq &\sum_{k\in
I_{h}}\int_{Y_{h}^{k}}\left| u_{h}\left( hx\right) \right| ^{p}dx  \nonumber
\\
&=&\sum_{k\in I_{h}}\left( \frac{1}{h}\right) ^{n}\int_{k+Y}\left|
u_{h}\left( x\right) \right| ^{p}dx  \label{ulikhet} \\
&\leq &K\int_{Y}\left| u_{h}\left( x\right) \right| ^{p}dx.  \nonumber
\end{eqnarray}
Now it follows that $u_{h}\left( h\cdot \right) $ is bounded in $%
L^{p}(\Omega )$ by taking into account that any weakly convergent sequence
is bounded. By H\"{o}lder's inequality we have that 
\begin{eqnarray*}
&&\left| \int_{\Omega }\left( u_{h}\left( hx\right) -\frac{1}{\left|
Y\right| }\int_{Y}u\left( x\right) dx\right) \phi (x)dx\right| \\
&\leq &\left( \int_{\Omega }\left| u_{h}\left( hx\right) -\frac{1}{\left|
Y\right| }\int_{Y}u\left( x\right) dx\right| ^{p}dx\right) ^{\tfrac{1}{p}%
}\left( \int_{\Omega }\left| \phi (x)-\phi ^{h}(x)\right| ^{q}dx\right) ^{%
\tfrac{1}{q}} \\
&&+\left| \int_{\Omega }\left( u_{h}\left( hx\right) -\frac{1}{\left|
Y\right| }\int_{Y}u_{h}(x)dx\right) \phi ^{h}(x)dx\right| \\
&&+\left| \frac{1}{\left| Y\right| }\int_{Y}\left( u_{h}(x)-u(x)\right) \phi
^{h}(x)\,dx\right| \left| \Omega \right| .
\end{eqnarray*}
This together with (\ref{lqconverg}), (\ref{large}) and (\ref{ulikhet})
implies that 
\[
\int_{\Omega }\left( u_{h}\left( hx\right) -\frac{1}{\left| Y\right| }%
\int_{Y}u\left( x\right) dx\right) \phi (x)dx\rightarrow 0, 
\]
as $h\rightarrow \infty $ for every $\phi \in C_{0}(\Omega )$. By using a
density argument it also holds for every $\phi \in L^{q}(\Omega )$ and we
are done.

Let us now turn to the case $p=1$. Let $u_{h}^{i}$ be defined as follows 
\[
u_{h}^{i}=\left\{ 
\begin{array}{l}
u_{h}\text{ \ \ if }u_{h}(x)<t_{1/i} \\ 
\\ 
0\;\;\;\text{ if }u_{h}(x)\geq t_{1/i}
\end{array}
\right. \text{\thinspace }. 
\]
According to Lemma \ref{lemdunford} there exists a constant $t_{1/i}>0$ for
each positive integer $i$ such that 
\begin{equation}
\int_{\Omega }\left| u_{h}^{i}-u_{h}\right| dx=\int_{\left\{ u_{h}\geq
t_{1/i}\right\} }\left| u_{h}\right| dx\leq \frac{1}{i},  \label{dag1}
\end{equation}
for all $h,i\in N$. By a diagonalization argument each subsequence of $(h)$
has a subsequence, denoted by $(h^{^{\prime }})$, such that $u_{h^{^{\prime
}}}^{i}$ converges weak* in $L^{\infty }(Y)$ to some function $u^{i}$ for
every $i$. It is easy to see that the proof for the case ($1<p\leq \infty $)
also holds if $(h)$ is replaced with $(h^{\prime })$, which implies that 
\begin{equation}
u_{h^{\prime }}^{i}(h^{\prime }\cdot )\rightarrow \frac{1}{\left| Y\right| }%
\int_{Y}u^{i}(x)dx  \label{h-merked}
\end{equation}
weak* in $L^{\infty }(\Omega )$ for every $i$. Let $v\in L^{\infty }(\Omega
) $. Then 
\begin{eqnarray}
&&\limsup_{h^{\prime }\rightarrow \infty }\left| \int_{\Omega }v(x)\left(
u_{h^{\prime }}\left( h^{\prime }x\right) dx-\left( \frac{1}{\left| Y\right| 
}\int_{Y}u(x)dx\right) \right) dx\right|  \label{limsup} \\
&=&\limsup_{i\rightarrow \infty }\limsup_{h^{\prime }\rightarrow \infty
}\left| \int_{\Omega }v(x)\left( u_{h^{\prime }}\left( h^{\prime }x\right)
dx-\left( \frac{1}{\left| Y\right| }\int_{Y}u(x)dx\right) \right) dx\right| 
\nonumber \\
&\leq &\limsup_{i\rightarrow \infty }\limsup_{h^{\prime }\rightarrow \infty
}\left| \int_{\Omega }v(x)\left( u_{h^{\prime }}\left( h^{\prime }x\right)
-u_{h^{\prime }}^{i}\left( h^{\prime }x\right) \right) dx\right|  \nonumber
\\
&&+\limsup_{i\rightarrow \infty }\limsup_{h^{\prime }\rightarrow \infty
}\left| \int_{\Omega }v(x)\left( u_{h^{\prime }}^{i}\left( h^{\prime
}x\right) dx-\left( \frac{1}{\left| Y\right| }\int_{Y}u(x)dx\right) \right)
dx\right| .  \nonumber
\end{eqnarray}
Both of the last terms are zero. For the first term this is seen by
replacing $u_{h}$ with $u_{h^{\prime }}^{i}-u_{h^{\prime }}$ in (\ref
{ulikhet}) and using (\ref{h-merked}) and (\ref{dag1}) to obtain that 
\begin{eqnarray*}
\left| \int_\Omega v(x)\left( u_{h^{\prime }}\left( h^{\prime }x\right)
-u_{h^{\prime }}^i\left( h^{\prime }x\right) \right) dx\right|  &\leq
&\left\| v(x)\right\| _\infty \int_\Omega \left| u_{h^{\prime }}^i\left(
h^{\prime }x\right) -u_{h^{\prime }}\left( h^{\prime }x\right) \right| dx \\
&\leq &\left\| v(x)\right\| _\infty K\int_Y\left| u_{h^{\prime }}^i\left(
x\right) -u_{h^{\prime }}\left( x\right) \right| dx \\
&\leq &\left\| v(x)\right\| _\infty \frac Ki.
\end{eqnarray*}
From this it is clear that 
\[
\limsup_{i\rightarrow \infty }\limsup_{h^{\prime }\rightarrow \infty }\left|
\int_{\Omega }v(x)\left( u_{h^{\prime }}\left( h^{\prime }x\right)
-u_{h^{\prime }}^{i}\left( h^{\prime }x\right) \right) dx\right| =0. 
\]
For the second term we use (\ref{h-merked}), the weak lower semicontinuity
of the $L^{1}(\Omega )$ norm and (\ref{dag1}) in order to obtain that 
\begin{eqnarray*}
&&\limsup_{h^{\prime }\rightarrow \infty }\left| \int_{\Omega }v(x)\left(
u_{h^{\prime }}^{i}\left( h^{\prime }x\right) dx-\left( \frac{1}{\left|
Y\right| }\int_{Y}u(x)dx\right) \right) dx\right| \\
&\leq &\frac{1}{\left| Y\right| }\int_{Y}\left| u^{i}(x)-u(x)\right|
dx\left| \int_{\Omega }v(x)dx\right| \\
&\leq &\frac{1}{\left| Y\right| }\liminf_{h^{\prime }\rightarrow \infty
}\int_{Y}\left| u_{h^{\prime }}^{i}\left( x\right) -u_{h^{\prime }}\left(
x\right) \right| dx\left| \int_{\Omega }v(x)dx\right| \\
&\leq &\frac{1}{\left| Y\right| i}\left| \int_{\Omega }v(x)dx\right|
\rightarrow 0
\end{eqnarray*}
as $i\rightarrow \infty $. Summing up from (\ref{limsup}) we have that any
subsequence of $\left( u_{h}(h\cdot )\right) $ contains a subsequence $%
\left( u_{h^{\prime }}(h^{^{\prime }}\cdot )\right) $ which converges weakly
to $\left| Y\right| ^{-1}\int_{Y}u(x)dx$ in $L^{1}(\Omega ).$ Thus this is
also true for the whole sequence $\left\{ u_{h}(h\cdot )\right\} $.

\section{Homogenization of some periodic boundary value problems}

Before we state the result of this section we introduce some definitions and
notations. Let $Y$ and $Z$ be a open bounded rectangles in $R^{n}$, $\left|
E\right| $ denotes the Lebesgue measure of the set $E\subset R^{n}$ and $%
(\cdot ,\cdot )$ is the Euclidean scalar product on $R^{n}$. Moreover, $c$
will be a constant that may differ from one place to an other and $h\in N$.
The function $\widetilde{\omega }:R\rightarrow R$ is an arbitrary function
which is continuous, increasing and $\widetilde{\omega }(0)=0$. By $%
W_{per}^{1,p}(Y)$ we denote the set of all functions $u\in W^{1,p}(Y)$ with
mean value zero which have the same trace on opposite faces of $Y$, $%
W_{per}^{1,p}(Z)$ is defined in the corresponding way. Every function $u\in
W_{per}^{1,p}(Y)$ can be extended by periodicity to a function in $%
W_{loc}^{1,p}(R^{n})$ (in this paper we will not make any distinction
between the original function and its extension). Let us fix a function $%
a:Y\times R^{n}\times R^{n}\rightarrow R^{n}$ which fulfills the conditions:

\begin{enumerate}
\item  $a(y,\cdot ,\xi )$ is $Z$-periodic and Lebesgue measurable for every $%
\xi \in R^n$ and every $y\in R^n,$

\item  There exists two constants $c_1,c_2>0$ and two constants $\alpha $
and $\beta $, with $0\leq \alpha \leq \min \left\{ 1,p-1\right\} $ and $\max
\left\{ p,2\right\} \leq \beta <\infty $ such that $a$ satisfies the
following boundedness, continuity and monotonicity assumptions: 
\begin{equation}
a(y,z,0)=0\;\;\text{for a.e. }y,z\in R^n,  \label{ayz0}
\end{equation}
\begin{equation}
\left| a(y,z,\xi _1)-a(y,z,\xi _2)\right| \leq c_1(1+\left| \xi _1\right|
+\left| \xi _2\right| )^{p-1-\alpha }\left| \xi _1-\xi _2\right| ^\alpha ,
\label{acont}
\end{equation}
\begin{equation}
(a(y,z,\xi _1)-a(y,z,\xi _2),\xi _1-\xi _2)\geq c_2(1+\left| \xi _1\right|
+\left| \xi _2\right| )^{p-\beta }\left| \xi _1-\xi _2\right| ^\beta
\label{amon}
\end{equation}

\item  $a$ is on the form $a(y,z,\xi )=\sum_{i=1}^N\chi _{\Omega
_i}(y)a_i(y,z,\xi )$ and satisfies a continuity condition of the form 
\begin{equation}
\left| a(y_1,z,\xi )-a(y_2,z,\xi )\right| ^q\leq \omega (\left|
y_1-y_2\right| )\left( 1+\left| \xi \right| ^p\right) ,\text{ \ \ }
\label{a1cont}
\end{equation}
for $y_1,y_2\in \Omega _i$, $i=1,\ldots ,N$,$\;$a.e. $z\in R^n$ and every $%
\xi \in R^n$, and where $\omega :R\rightarrow R$ is continuous, increasing
and $\omega (0)=0$.
\end{enumerate}

By (\ref{ayz0}), (\ref{acont}), and (\ref{amon}) it follows that 
\begin{equation}
\left| a(y,z,\xi )\right| \leq c\left( 1+\left| \xi \right| ^{p-1}\right) ,
\label{a1}
\end{equation}
\begin{equation}
\left| \xi \right| ^{p}\leq c\left( 1+(a(y,z,\xi ),\xi )\right) ,  \label{a2}
\end{equation}
hold for $y\in R^{n},\;$a.e. $z\in R^{n}$ and every $\xi \in R^{n}$.

We are now in the position to state the result in this section.

\begin{theorem}
\label{lemaux2}Let $a$ satisfy (\ref{ayz0}), (\ref{acont}), (\ref{amon}) and
(\ref{a1cont}). Moreover, let $(u_h^\xi )$ be the solutions of 
\begin{equation}
\left\{ 
\begin{array}{l}
\int_Y(a(x,hx,\xi +Du_h^\xi ),D\phi )\,dx=0\text{ for every }\phi \in
W_{per}^{1,p}(Y), \\ 
\\ 
u_h^\xi \in W_{per}^{1,p}(Y).
\end{array}
\right.  \label{kalle}
\end{equation}
Then 
\begin{eqnarray*}
u_h^\xi &\rightarrow &u^\xi \text{ weakly in }W_{per}^{1,p}(Y), \\
a(x,hx,\xi +Du_h^\xi ) &\rightarrow &b(x,\xi +Du^\xi )\text{ weakly in }%
L^q(Y,R^n),
\end{eqnarray*}
as $h\rightarrow \infty $, where $u^\xi $ is the unique solution of 
\[
\left\{ 
\begin{array}{l}
\int_Y(b(x,\xi +Du^\xi ),D\phi )\,dx=0\text{ for every }\phi \in
W_{per}^{1,p}(Y), \\ 
\\ 
u^\xi \in W_{per}^{1,p}(Y).
\end{array}
\right. 
\]
The operator $b:Y\times R^n\rightarrow R^n$ is defined as 
\[
b(y,\tau )=\frac 1{\left| Z\right| }\int_Za(y,z,\tau +Dv^{\tau ,y}(z))\,dz, 
\]
where $v^{\tau ,y}$ is the unique solution of the cell-problem 
\begin{equation}
\left\{ 
\begin{array}{l}
\int_Z(a(y,z,\tau +Dv^{\tau ,y}(z)),D\phi )\,dz=0\text{ for every }\phi \in
W_{per}^{1,p}(Z), \\ 
\\ 
v^{\tau ,y}\in W_{per}^{1,p}(Z).
\end{array}
\right.  \label{localeq1}
\end{equation}
\end{theorem}

%TCIMACRO{\TeXButton{Proof}{\proof}}
%BeginExpansion
\proof%
%EndExpansion
We divide the proof into several steps.

\textbf{Step 1.} Let $\left\{ \Omega _i^k\subset \Omega :i\in I_k\right\} $
denote a family of disjoint open sets with diameter less than $\frac 1k$
such that $\left| \Omega \backslash \cup _{i\in I_k}\Omega _i^k\right| =0$
and $\left| \partial \Omega _i^k\right| =0$. We define the function $a^k$ as 
\[
a^k(y,z,\xi )=\sum_{i\in I_k}\chi _{\Omega _i^k}(y)a(y_i^k,z,\xi ), 
\]
where $y_i^k\in \Omega _i^k$. Consider the auxillary periodic boundary value
problems (transmission problems) 
\begin{equation}
\left\{ 
\begin{array}{l}
\int_Y(a^k(x,hx,\xi +Du_h^{k,\xi }),D\phi )\,dx=0\text{\thinspace \thinspace
\thinspace \thinspace \thinspace \thinspace for every }\phi \in
W_{per}^{1,p}(Y), \\ 
\\ 
u_h^{k,\xi }\in W_{per}^{1,p}(Y).
\end{array}
\right.  \label{uhk}
\end{equation}
Then we have that 
\[
u_h^{k,\xi }\rightarrow u^{k,\xi }\text{ weakly in }W_{per}^{1,p}(Y), 
\]
\[
a^k(x,hx,\xi +Du_h^{k,\xi })\rightarrow b^k(x,\xi +Du^{k,\xi })\text{ weakly
in }L^q(Y,R^n), 
\]
where $u^{k,\xi }$ is the unique solution of the homogenized problem 
\begin{equation}
\left\{ 
\begin{array}{l}
\int_Y(b^k(x,\xi +Du^{k,\xi }),D\phi )\,dx=0\text{ for every }\phi \in
W_{per}^{1,p}(Y), \\ 
\\ 
u^{k,\xi }\in W_{per}^{1,p}(Y).
\end{array}
\right.  \label{homprob}
\end{equation}
The operator $b^k:Y\times R^n\rightarrow R^n$ is defined a.e. as 
\[
b^k(y,\tau )=\sum_{i=I_k}\chi _{\Omega _i^k}(y)\int_Za(y_i^k,z,\tau
+Dv^{\tau ,y_i^k}(z))\,dz=\sum_{i=I_k}\chi _{\Omega _i^k}(y)b(y_i^k,\tau ), 
\]
where $v^{\tau ,y_i^k}$ is the unique solution of the cell problem 
\begin{equation}
\left\{ 
\begin{array}{l}
\int_Z(a(y_i^k,z,\tau +Dv^{\tau ,y_i^k}(z)),D\phi (z))\,dz=0\text{ for every 
}\phi \in W_{per}^{1,p}(Z), \\ 
\\ 
v^{\tau ,y_i^k}\in W_{per}^{1,p}(Z).
\end{array}
\right.  \label{cellproblem3}
\end{equation}
The proof of these convergence results follows by suitable modifications of
well-known homogenization techniques. Indeed, let $\phi =u_h^{k,\xi }$ in (%
\ref{kalle}) then it follows by (\ref{a2}), (\ref{ayz0}), (\ref{acont}) and
H\"{o}lder's inequality that 
\begin{eqnarray}
\int_Y\left| \xi +Du_h^{k,\xi }\right| ^p\,dx &\leq &c\int_Y1+(a^k(x,hx,\xi
+Du_h^{k,\xi }),\xi +Du_h^{k,\xi })\,dx  \nonumber \\
&=&c\int_Y1+(a^k(x,hx,\xi +Du_h^{k,\xi }),\xi )\,dx  \nonumber \\
&\leq &c\left( 1+\int_Y\left( 1+\left| \xi +Du_h^{k,\xi }\right| \right)
^{p-1}\,dx\right)  \label{uhxibound1} \\
&\leq &c\left( 1+\left( \int_Y\left( \left| \xi +Du_h^{k,\xi }\right|
\right) ^p\,dx\right) ^{^{\tfrac 1q}}\right) .  \nonumber
\end{eqnarray}
If $\left( \int_Y\left( \left| \xi +Du_h^{k,\xi }\right| \right)
^p\,dx\right) ^{1/q}\leq 1$ it is clear that the sequence of solutions $%
(u_h^{k,\xi })$ is bounded in $L^p(Y,R^n)$, so let us assume that $\left(
\int_Y\left( \left| \xi +Du_h^{k,\xi }\right| \right) ^p\,dx\right)
^{1/q}\geq 1$, then (\ref{uhxibound1}) implies that 
\[
\int_Y\left| \xi +Du_h^{k,\xi }\right| ^p\,dx\leq c 
\]
which means that $(u_h^{k,\xi })$ is bounded in . Since $\left\| D\cdot
\right\| _{L^p(Y,R^n)}$ is an equivalent norm on $W_{per}^{1,p}(Y)$ it
follows that there exists a constant $c>0$ independent of $h$ such that 
\[
\left\| u_h^{k,\xi }\right\| _{W_{per}^{1,p}(Y)}\leq c. 
\]
From the reflexivity of $W_{per}^{1,p}(Y)$ there exists a subsequence, still
denoted by $(u_h^{k,\xi })$ such that 
\[
u_h^{k,\xi }\rightarrow u_{*}^{k,\xi }\text{ weakly in }W_{per}^{1,p}(Y). 
\]
Let us now define 
\[
\eta _h^{i,k,\xi }=a(x_i^k,hx,\xi +Du_h^{k,\xi }),\;\;\;i\in I_k 
\]
By (\ref{acont}), (\ref{ayz0}), H\"{o}lder's inequality and (\ref{uhxibound1}%
) we have that $\eta _h^{i,k,\xi }$ is bounded in $L^q(\Omega _i^k,R^n)$.
Indeed 
\begin{eqnarray*}
\int_{\Omega _i^k}\left| \eta _h^{i,k,\xi }\right| ^q\,dx &=&\int_{\Omega
_i^k}\left| a(x_i^k,hx,\xi +Du_h^{k,\xi })\right| ^q\,dx \\
&\leq &c\int_{\Omega _i^k}\left( 1+\left| \xi +Du_h^{k,\xi })\right| \right)
^{q(p-1-\alpha )}\left| \xi +Du_h^{k,\xi })\right| ^{q\alpha }\,dx \\
&\leq &c\int_{\Omega _i^k}1+\left| \xi +Du_h^{k,\xi })\right| ^p\,dx\leq c
\end{eqnarray*}
where $c$ is a constant independent of $h$. This means that there exists a
subsequence, still denoted by $(\eta _h^{i,k,\xi }),$ and a $\eta
_{*}^{i,k,\xi }\in L^q(\Omega _i^k,R^n)$ such that 
\[
\eta _h^{i,k,\xi }\rightarrow \eta _{*}^{i,k,\xi }\text{ weakly in }%
L^q(\Omega _i^k,R^n). 
\]
From our original problem (\ref{kalle}) we have that 
\[
\left\{ 
\begin{array}{l}
\sum_{i\in I_k}\int_{\Omega _i^k}(a(x_i^k,hx,\xi +Du_h^{k,\xi }),D\phi
)\,dx=0\text{ for every }\phi \in W_{per}^{1,p}(Y), \\ 
\\ 
u_h^{k,\xi }\in W_{per}^{1,p}(Y).
\end{array}
\right. 
\]
In the limit we get 
\[
\sum_{i\in I_k}\int_{\Omega _i^k}(\eta _{*}^{i,k,\xi },D\phi )\,dx=0\text{
for every }\phi \in W_{per}^{1,p}(Y). 
\]
Especially this means that 
\[
\int_{\Omega _i^k}(\eta _{*}^{i,k,\xi },D\phi )\,dx=0\text{ for every }\phi
\in C_0^\infty (\Omega _i^k)\text{, \thinspace \thinspace \thinspace
\thinspace \thinspace \thinspace \thinspace }i\in I_k. 
\]
If we now could show that 
\begin{equation}
\eta _{*}^{i,k,\xi }=b(x_i^k,\xi +Du_{*}^{k,\xi })\text{ for a.e }x\in
\Omega _i^k,  \label{eq15}
\end{equation}
then it follows by the uniqueness of the homogenized problem (\ref{homprob})
that $u_{*}^{k,\xi }=u^{k,\xi }$. To this aim we define the function 
\[
w_h^{\tau ,x_i^k}(x)=(\tau ,x)+\frac 1hv^{\tau ,x_i^k}(hx), 
\]
where $v^{\tau ,x_i^k}$ is defined as in (\ref{cellproblem3}). By
periodicity we have that 
\begin{eqnarray*}
w_h^{\tau ,x_i^k} &\rightarrow &(\tau ,\cdot )\text{ weakly in }%
W^{1,p}(\Omega _i^k), \\
Dw_h^{\tau ,x_i^k} &\rightarrow &\tau \text{ weakly in }L^p(\Omega _i^k,R^n),
\\
a(x_i^k,hx,Dw_h^{\tau ,x_i^k}) &\rightarrow &b(x_i^k,\tau )\text{ weakly in }%
L^q(\Omega _i,R^n).
\end{eqnarray*}
By the monotonicity of $a_i$ we have for a fix $\tau $ that 
\[
\int_{\Omega _i}(a(x_i^k,hx,\xi +Du_h^{k,\xi })-a(x_i^k,hx,Dw_h^{\tau
,x_i^k}),\xi +Du_h^{k,\xi }-Dw_h^{\tau ,x_i^k})\phi \,dx\geq 0, 
\]
for every $\phi \in C_0^\infty (\Omega _i),\phi \geq 0.$ By density we
obtain that 
\[
(\eta _{*}^{i,k,\xi }(x)-b(x_i^k,\tau ),\xi +Du_{*}^{k,\xi }(x)-\tau )\geq 0%
\text{ for a.e. }x\in \Omega _i^k\text{ and for every }\tau \in R^n. 
\]
Since $b^k$ is monotone and continuous, see Proposition \ref{lemmab1prop},
we have that $b^k$ is maximal monotone and the crucial relation (\ref{eq15})
follows. We have now proved step 1 up to a subsequence of $(u_h^{k,\xi })$.
By the uniqueness of the solution of the homogenized equation (\ref{homprob}%
) it follows that it is true for the whole sequence.

\textbf{Step 2.} Let us now prove that $u_{h}^{\xi }\rightarrow u^{\xi }$
weakly in $W_{per}^{1,p}(Y)$. Let $g\in \left( W_{per}^{1,p}(Y)\right) ^{*},$
then 
\begin{eqnarray*}
\lim_{h\rightarrow \infty }\left\langle g,u_{h}^{\xi }-u^{\xi }\right\rangle
&=&\lim_{k\rightarrow \infty }\lim_{h\rightarrow \infty }\left\langle
g,u_{h}^{\xi }-u^{\xi }\right\rangle \\
&\leq &\lim_{k\rightarrow \infty }\lim_{h\rightarrow \infty }\left\|
g\right\| _{\left( W_{per}^{1,p}(Y)\right) ^{*}}\left\| u_{h}^{\xi
}-u_{h}^{k,\xi }\right\| _{W_{per}^{1,p}(Y)} \\
&&+\lim_{k\rightarrow \infty }\lim_{h\rightarrow \infty }\left\langle
g,u_{h}^{k,\xi }-u^{k,\xi }\right\rangle \\
&&+\lim_{k\rightarrow \infty }\left\| g\right\| _{\left(
W_{per}^{1,p}(Y)\right) ^{*}}\left\| u^{k,\xi }-u^{\xi }\right\|
_{W_{per}^{1,p}(Y)}.
\end{eqnarray*}
It is enough to prove that all three terms on the right hand side are zero.

\textbf{Term 1. }Let us prove that 
\begin{equation}
\lim_{k\rightarrow \infty }\lim_{h\rightarrow \infty }\left\| u_h^\xi
-u_h^{k,\xi }\right\| _{W_{per}^{1,p}(Y)}=0.  \label{step1a}
\end{equation}
By definition 
\begin{eqnarray*}
\int_Y(a^k(x,hx,\xi +Du_h^{k,\xi }),D\phi )\,dx &=&0\text{ for every }\phi
\in W_{per}^{1,p}(Y), \\
\int_Y(a(x,hx,\xi +Du_h^\xi ),D\phi )\,dx &=&0\text{ for every }\phi \in
W_{per}^{1,p}(Y).
\end{eqnarray*}
This implies that we for $\phi =u_h^{k,\xi }-u_h^\xi $ have 
\begin{eqnarray*}
&&\int_Z(a^k(x,hx,\xi +Du_h^{k,\xi })-a^k(x,hx,\xi +Du_h^\xi ),Du_h^{k,\xi
}-Du_h^\xi )\,dx \\
&=&\int_Z(a(x,hx,\xi +Du_h^\xi )-a^k(x,hx,\xi +Du_h^\xi ),Du_h^{k,\xi
}-Du_h^\xi )\,dx.
\end{eqnarray*}
By using (\ref{amon}), and H\"{o}lder's reversed inequality on the left hand
side and H\"{o}lder inequality (\ref{a1cont}) and the fact that $(u_h^\xi )$
and $(u_h^{k,\xi })$ is bounded in $W_{per}^{1,p}(Y)$ on the right hand side
we obtain that 
\begin{eqnarray*}
&&c_2\left( \int_Y\left| Du_h^{k,\xi }-Du_h^\xi \right| ^p\,dx\right)
^{\tfrac \beta p} \\
&&\times \left( \int_Y\left( 1+\left| \xi +Du_h^{k,\xi }\right| +\left| \xi
+Du_h^\xi \right| \right) ^p\,dx\right) ^{\tfrac p{p-\beta }} \\
&\leq &c_2\int_Y\left( 1+\left| \xi +Du_h^{k,\xi }\right| +\left| \xi
+Du_h^\xi \right| \right) ^{p-\beta }\left| Du_h^{k,\xi }-Du_h^\xi \right|
^\beta \,dx \\
&\leq &\left( \int_Y\left| (a(x,\frac x{\varepsilon _h},\xi +Du_h^\xi
)-a^k(x,\frac x{\varepsilon _h},\xi +Du_h^\xi )\right| ^q\,dx\right)
^{\tfrac 1q} \\
&&\times \left( \int_Y\left| Du_h^{k,\xi }-Du_h^\xi \right| ^p\,dx\right)
^{\tfrac 1p} \\
&\leq &\widetilde{\omega }(\frac 1k)\left( \int_Y1+\left| \xi +Du_h^\xi
)\right| ^p\,dx\right) ^{\tfrac 1q}\left( \int_Y\left| Du_h^{k,\xi
}-Du_h^\xi \right| ^p\,dx\right) ^{\tfrac 1p} \\
&\leq &\widetilde{\omega }(\frac 1k)\left( \int_Y\left| Du_h^{k,\xi
}-Du_h^\xi \right| ^p\,dx\right) ^{\tfrac 1p}.
\end{eqnarray*}
Since $\left\| D\cdot \right\| _{L^p(Y,R^n)}$ is an equivalent norm on $%
W_{per}^{1,p}(Y)$ this implies that 
\begin{equation}
\left\| u_h^{k,\xi }-u_h^\xi \right\| _{W_{per}^{1,p}(Y)}\leq \widetilde{%
\omega }(\frac 1k)\rightarrow 0  \label{bound}
\end{equation}
as $k\rightarrow \infty $ uniformly in $h$. This means that we can change
the order in the limit process in (\ref{step1a}) and (\ref{step1a}) follows
by taking (\ref{bound}) into account.

\textbf{Term 2.} We observe that 
\[
\lim_{k\rightarrow \infty }\lim_{h\rightarrow \infty }\left\langle
g,u_h^{k,\xi }-u^{k,\xi }\right\rangle =0, 
\]
as a direct consequence of Step 1.

\textbf{Term 3.} Let us prove that 
\begin{equation}
\lim_{k\rightarrow \infty }\left\| u^{k,\xi }-u^{\xi }\right\|
_{W_{per}^{1,p}(Y)}=0.  \label{step3}
\end{equation}
By definition we have that 
\begin{eqnarray*}
\int_{Y}(b^{k}(x,\xi +Du^{k,\xi }),D\phi )\,dx &=&0\text{ for every }\phi
\in W_{per}^{1,p}(Y), \\
\int_{Y}(b(x,\xi +Du^{\xi }),D\phi )\,dx &=&0\text{ for every }\phi \in
W_{per}^{1,p}(Y).
\end{eqnarray*}
Thus 
\begin{eqnarray*}
&&\int_{Y}(b^{k}(x,\xi +Du^{k,\xi })-b^{k}(x,\xi +Du^{\xi }),D\phi )\,dx \\
&=&\int_{Y}(b(x,\xi +Du^{\xi })-b^{k}(x,\xi +Du^{\xi }),D\phi )\,dx,
\end{eqnarray*}
for every $\phi \in W_{per}^{1,p}(Y)$. Choose $\phi =u^{k,\xi }-u^{\xi }$
and take the strict monotonicity of $b^{k}$, see (\ref{bmon2}), into account
on the left hand side and apply the H\"{o}lder inequality and (\ref{bcont2})
on the right hand side to obtain 
\begin{eqnarray*}
&&c\left( \int_{Y}\left| Du^{k,\xi }-Du^{\xi }\right| ^{p}\,dx\right) ^{%
\tfrac{\beta }{p}}\times \\
&&\left( \int_{Y}\left( 1+\left| \xi +Du^{k,\xi }\right| +\left| \xi
+Du^{\xi }\right| \right) ^{p}\,dx\right) ^{\tfrac{p-\beta }{p}} \\
&\leq &c\int_{Y}\left( 1+\left| \xi +Du^{k,\xi }\right| +\left| \xi +Du^{\xi
}\right| \right) ^{p-\beta }\left| Du^{k,\xi }-Du^{\xi }\right| ^{\beta }\,dx
\\
&\leq &\left( \int_{Y}\left| b(x,\xi +Du^{\xi })-b^{k}(x,\xi +Du^{\xi
})\right| ^{q}\,dx\right) ^{\tfrac{1}{q}}\times \\
&&\left( \int_{Y}\left| Du^{k,\xi }-Du^{\xi }\right| ^{p}\,dx\right) ^{%
\tfrac{1}{p}} \\
&\leq &\widetilde{\omega }(\frac{1}{k})\left( \int_{Y}1+\left| \xi +Du^{\xi
}\right| ^{p}\,dx\right) ^{\tfrac{1}{q}}\left( \int_{Y}\left| Du^{k,\xi
}-Du^{\xi })\right| ^{p}\,dx\right) ^{\tfrac{1}{p}}.
\end{eqnarray*}
By using the fact that $u^{\xi }$ and $u^{k,\xi }$ are bounded in $%
W_{per}^{1,p}(Y)$ it follows that 
\begin{equation}
\left\| Du^{k,\xi }-Du^{\xi }\right\| _{L^{p}(Y,R^{n})}\leq \widetilde{%
\omega }(\frac{1}{k}),  \label{bound2}
\end{equation}
and the result follows by noting that $\left\| D\cdot \right\|
_{L^{p}(Y,R^{n})}$ is an equivalent norm on $W_{per}^{1,p}(Y)$.

\textbf{Step 3.} Next we prove that $a(x,hx,\xi +Du_{h}^{\xi })\rightarrow
b(x,\xi +Du^{\xi })$ weakly in $L^{q}(Y,R^{n})$. In fact if $g\in
(L^{q}(Y,R^{n}))^{*},$ then 
\begin{eqnarray*}
&&\lim_{h\rightarrow \infty }\left\langle g,a(x,hx,\xi +Du_{h}^{\xi
})-b(x,\xi +Du^{\xi })\right\rangle \\
&=&\lim_{k\rightarrow \infty }\lim_{h\rightarrow \infty }\left\langle
g,a(x,hx,\xi +Du_{h}^{\xi })-b(x,\xi +Du^{\xi })\right\rangle \\
&\leq &\lim_{k\rightarrow \infty }\lim_{h\rightarrow \infty }\left\|
g\right\| \left\| a(x,hx,\xi +Du_{h}^{\xi })-a^{k}(x,hx,\xi +Du_{h}^{k,\xi
})\right\| _{L^{q}(Y,R^{n})} \\
&&+\lim_{k\rightarrow \infty }\lim_{h\rightarrow \infty }\left\langle
g,a^{k}(x,hx,\xi +Du_{h}^{k,\xi })-b^{k}(x,\xi +Du^{k,\xi }))\right\rangle \\
&&+\lim_{k\rightarrow \infty }\left\| g\right\| \left\| b^{k}(x,\xi
+Du^{k,\xi })-b(x,\xi +Du^{\xi })\right\| _{L^{q}(Y,R^{n})}.
\end{eqnarray*}
It is sufficient to prove that all three terms on the right hand side are
zero.

\textbf{Term 1.} Let us show that 
\begin{equation}
\lim_{k\rightarrow \infty }\lim_{h\rightarrow \infty }\left\| a(x,hx,\xi
+Du_{h}^{\xi })-a^{k}(x,hx,\xi +Du_{h}^{k,\xi })\right\| _{L^{q}(Y,R^{n})}=0.
\label{step12}
\end{equation}
By using elementary estimates we find that 
\begin{eqnarray*}
&&\int_{Y}\left| a^{k}(x,hx,\xi +Du_{h}^{k,\xi })-a(x,hx,\xi +Du_{h}^{\xi
})\right| ^{q}\,dx \\
&\leq &c\int_{Y}\left| a^{k}(x,hx,\xi +Du_{h}^{k,\xi })-a^{k}(x,hx,\xi
+Du_{h}^{\xi })\right| ^{q}\,dx \\
&&+c\int_{Y}\left| a^{k}(x,hx,\xi +Du_{h}^{\xi })-a(x,hx,\xi +Du_{h}^{\xi
})\right| ^{q}\,dx.
\end{eqnarray*}
Hence, by applying the continuity conditions (\ref{acont}) and H\"{o}lder
inequality to the first term and (\ref{a1cont}) to the second term we obtain
that 
\begin{eqnarray*}
&&\int_{Y}\left| a^{k}(x,hx,\xi +Du_{h}^{k,\xi })-a(x,hx,\xi +Du_{h}^{\xi
})\right| ^{q}\,dx \\
&\leq &c\left( \int_{Y}\left( 1+\left| \xi +Du_{h}^{k,\xi }\right| +\left|
\xi +Du_{h}^{\xi }\right| \right) ^{p}\,dx\right) ^{\tfrac{p-1-\alpha }{p-1}}
\\
&&\times \left( \int_{Y}\left| Du_{h}^{k,\xi }-Du_{h}^{\xi }\right|
^{p}\,dx\right) ^{\tfrac{\alpha }{p-1}}+\widetilde{\omega }(\frac{1}{k}%
)\int_{Y}1+\left| \xi +Du_{h}^{\xi }\right| ^{p}\,dx.
\end{eqnarray*}
By using the fact that $u_{h}^{k,\xi }$ and $u_{h}^{\xi }$ are bounded in $%
W_{per}^{1,p}(Y)$ and (\ref{bound}) it follows that 
\begin{equation}
\left\| a(x,hx,\xi +Du_{h}^{\xi })-a^{k}(x,hx,\xi +Du_{h}^{k,\xi })\right\|
_{L^{q}(Y,R^{n})}\leq \widetilde{\omega }(\frac{1}{k})\rightarrow 0.
\label{step121}
\end{equation}
as $k\rightarrow \infty $ uniformly in $h.$ This implies that we may change
the order in the limit process in (\ref{step12}) and we obtain (\ref{step12}%
) by taking (\ref{step121}) into account.

\textbf{Term 2.} We observe that 
\[
\lim_{k\rightarrow \infty }\lim_{h\rightarrow \infty }\left\langle
g,a^k(x,hx,\xi +Du_h^{k,\xi })-b^k(x,\xi +Du^{k,\xi })\right\rangle =0, 
\]
as a direct consequence of Step 1.

\textbf{Term 3}. Let us show that 
\begin{equation}
\lim_{k\rightarrow \infty }\left\| b^{k}(x,\xi +Du^{k,\xi })-b(x,\xi
+Du^{\xi })\right\| _{L^{q}(Y,R^{n})}=0.  \label{step13}
\end{equation}
We have that 
\begin{eqnarray*}
&&\int_{Y}\left| b^{k}(x,\xi +Du_{*}^{k,\xi })-b_{1}(x,\xi +Du^{\xi
})\right| ^{q}\,dx \\
&\leq &c\int_{Y}\left| b^{k}(x,\xi +Du_{*}^{k,\xi })-b^{k}(x,\xi +Du^{\xi
})\right| ^{q}\,dx \\
&&+c\int_{Y}\left| b^{k}(x,\xi +Du^{\xi })-b_{1}(x,\xi +Du^{\xi })\right|
^{q}\,dx.
\end{eqnarray*}
By applying the continuity condition (\ref{bcont12}) and H\"{o}lders's
inequality to the first term and the continuity condition (\ref{bcont2}) to
the second term we see that 
\begin{eqnarray*}
&&\int_{Y}\left| b^{k}(x,\xi +Du^{k,\xi })-b(x,\xi +Du^{\xi })\right|
^{q}\,dx \\
&\leq &c\left( \int_{Y}\left( 1+\left| \xi +Du^{k,\xi }\right| +\left| \xi
+Du^{\xi }\right| \right) ^{p}\,dx\right) ^{\tfrac{p-1-\gamma }{p-1}} \\
&&\times \left( \int_{Y}\left| Du^{k,\xi }-Du^{\xi }\right| ^{p}\,dx\right)
^{\tfrac{\gamma }{p-1}}+\widetilde{\omega }(\frac{1}{k})\int_{Y}\left|
Du^{\xi }\right| ^{p}\,dx.
\end{eqnarray*}
By using the fact that $u^{k,\xi }$ and $u^{\xi }$ are bounded in $%
W_{per}^{1,p}(Y)$ and (\ref{bound2}) it follows that 
\[
\left\| b^{k}(x,\xi +Du^{k,\xi })-b(x,\xi +Du^{\xi })\right\|
_{L^{q}(Y,R^{n})}\leq \widetilde{\omega }(\frac{1}{k})\rightarrow 0 
\]
and we are done.

%TCIMACRO{\TeXButton{hfill}{\hfill}}
%BeginExpansion
\hfill%
%EndExpansion
{}$\square $

We remark that we have only considered the case when $a$ satisfies (\ref
{a1cont}) over the whole $Y$ the piecewise case follows by using the
technique used in step 1.

\begin{proposition}
\label{lemmab1prop}Let $b$ be the homogenized operator defined in Theorem 
\ref{lemaux2}. Then
\end{proposition}

\begin{itemize}
\item[(i)]  $b(\cdot ,\xi )$ satisfies the continuity condition 
\begin{equation}
\left| b(y_1,\xi )-b(y_2,\xi )\right| ^q\leq \widetilde{\omega }(\left|
y_1-y_2\right| )\left( 1+\left| \xi \right| ^p\right) .  \label{bcont2}
\end{equation}

\item[(ii)]  $b(x,\cdot )$ is strictly monotone, more precisely 
\begin{equation}
(b_1(y,\xi _1)-b_1(y,\xi _2),\xi _1-\xi _2)\geq c\left( 1+\left| \xi
_1\right| +\left| \xi _2\right| \right) ^{p-\beta }\left| \ \xi _1-\xi
_2\right| ^\beta ,  \label{bmon2}
\end{equation}
$\xi _1,\xi _2\in R^n.$

\item[(iii)]  $b(x,\cdot )$ is Lipschitz continuous, more precisely 
\begin{equation}
\left| b(x,\xi _1)-b(x,\xi _2)\right| \leq c\left( 1+\left| \xi _1\right|
+\left| \xi _2\right| \right) ^{p-1-\gamma }\left| \ \xi _1-\xi _2\right|
^\gamma ,  \label{bcont12}
\end{equation}
for every $\xi _1,\xi _2\in R^n$, where $\gamma =\alpha /(\beta -\alpha )$.

\item[(iv)]  
\begin{equation}
b(x,0)=0\text{ for }x\in Z.  \label{b2x0}
\end{equation}
\end{itemize}

%TCIMACRO{\TeXButton{Proof}{\proof }}
%BeginExpansion
\proof %
%EndExpansion
(i): By the definition of $b$ and Jensen's inequality we have that 
\begin{eqnarray*}
&&\left| b(y_1,\tau )-b(y_2,\tau )\right| ^q \\
&=&\left| \int_Za(y_1,z,\tau +Dv^{\tau ,y_1}(z))\,-a(y_2,z,\tau +Dv^{\tau
,y_2}(z))\,dz\right| ^q \\
&\leq &c\int_Z\left| a(y_1,z,\tau +Dv^{\tau ,y_1}(z))\,-a(y_2,z,\tau
+Dv^{\tau ,y_1}(z))\right| ^q\,dz \\
&&+c\int_Z\left| a(y_2,z,\tau +Dv^{\tau ,y_1}(z))-a(y_2,z,\tau +Dv^{\tau
,y_2}(z))\right| ^q\,dz.
\end{eqnarray*}
By applying (\ref{a1cont}) to the first term and (\ref{acont}) in
combination with H\"{o}lder's inequality to the second term we obtain that 
\begin{eqnarray}
&&\left| b(y_1,\tau )-b(y_2,\tau )\right| ^q\leq \widetilde{\omega }(\left|
y_1-y_2\right| )\int_Z1+\left| \tau +Dv^{\tau ,y_1}\right| ^p\,dz  \nonumber
\\
+ &&c\left( \int_Z\left( 1+\left| \tau +Dv^{\tau ,y_1}\right| +\left| \tau
+Dv^{\tau ,y_2}\right| \right) ^p\,dz\right) ^{\tfrac{p-1-\alpha }{p-1}} 
\nonumber \\
&&\times \left( \int_Z\left| Dv^{\tau ,y_1}-Dv^{\tau ,y_2}\right|
^p\,dz\right) ^{\tfrac \alpha {p-1}}  \label{bcont1}
\end{eqnarray}
Let us now study the two terms in (\ref{bcont1}) separately. The first term:
(\ref{a2}), (\ref{localeq1}) and (\ref{a1}) yields

\begin{eqnarray*}
\int_Z\left| \tau +Dv^{\tau ,y_1}\right| ^p\,dz &\leq &c\int_Z1+(a(y,z,\tau
+Dv^{\tau ,y_1}),\tau +Dv^{\tau ,y_1})\,dz \\
&=&c\int_Z1+(a(y,z,\tau +Dv^{\tau ,y_1}),\tau )\,dz \\
&\leq &c\int_Z1+c\left( 1+\left| \tau +Dv^{\tau ,y_1}\right| ^{p-1}\right)
\left| \tau \right| \,dz.
\end{eqnarray*}
By using the Young inequality we obtain that 
\begin{equation}
\int_Z\left| \tau +Dv^{\tau ,y_1}\right| ^p\,dz\leq c\left( 1+\left| \tau
\right| ^p\right) .  \label{bconta}
\end{equation}
Let us now study the second term in (\ref{bcont1}): By definition we have
that 
\begin{eqnarray*}
\int_Z(a(y_1,z,\tau +Dv^{\tau ,y_1}),D\phi )\,dz &=&0\text{ for every }\phi
\in W_{per}^{1,p}(Z), \\
\int_Z(a(y_2,z,\tau +Dv^{\tau ,y_2}),D\phi )\,dz &=&0\text{ for every }\phi
\in W_{per}^{1,p}(Z).
\end{eqnarray*}
This implies that 
\begin{eqnarray*}
&&\int_Z(a(y_1,z,\tau +Dv^{\tau ,y_1})-a(y_1,z,\tau +Dv^{\tau ,y_2}),D\phi
)\,dz \\
&=&\int_Z(a(y_2,z,\tau +Dv^{\tau ,y_2})-a(y_1,z,\tau +Dv^{\tau ,y_2}),D\phi
)\,dz,
\end{eqnarray*}
for every $\phi \in W_{per}^{1,p}(Z)$. In particular, for $\phi =v^{\tau
,y_1}-v^{\tau ,y_2}$, we have that 
\begin{eqnarray*}
&&\int_Z(a(y_1,z,\tau +Dv^{\tau ,y_1})-a(y_1,z,\tau +Dv^{\tau
,y_2}),Dv^{\tau ,y_1}-Dv^{\tau ,y_2})\,dz \\
&=&\int_Z(a(y_2,z,\tau +Dv^{\tau ,y_2})-a(y_1,z,\tau +Dv^{\tau
,y_2}),Dv^{\tau ,y_1}-Dv^{\tau ,y_2})\,dz.
\end{eqnarray*}
By applying the reversed H\"{o}lder inequality and (\ref{amon}) on the left
hand side and Schwarz's and H\"{o}lder's inequalities on the right hand side
it follows that 
\begin{eqnarray*}
&&c\left( \int_Z\left| Dv^{\tau ,y_1}-Dv^{\tau ,y_2}\right| ^p\;dz\right)
^{\tfrac \beta p}\times \\
&&\left( \int_Z\left( 1+\left| \tau +Dv^{\tau ,y_1}\right| +\left| \tau
+Dv^{\tau ,y_2}\right| \right) ^p\,dz\right) ^{\tfrac{p-\beta }p} \\
&\leq &c\int_Z\left( 1+\left| \tau +Dv^{\tau ,y_1}\right| +\left| \tau
+Dv^{\tau ,y_2}\right| \right) ^{p-\beta }\left| Dv^{\tau ,y_1}-Dv^{\tau
,y_2}\right| ^\beta \,dz \\
&\leq &\left( \int_Z\left| a(y_2,z,\tau +Dv^{\tau ,y_2})-a(y_1,z,\tau
+Dv^{\tau ,y_2})\right| ^q\,dz\right) ^{\tfrac 1q} \\
&&\times \left( \int_Z\left| Dv^{\tau ,y_1}-Dv^{\tau ,y_2}\right|
^p\,dz\right) ^{\tfrac 1p},
\end{eqnarray*}
which means that 
\begin{eqnarray}
&&\left( \int_Z\left| Dv^{\tau ,y_1}-Dv^{\tau ,y_2}\right| ^p\;dz\right)
^{\tfrac \alpha {p-1}}  \nonumber \\
&\leq &c\left( \int_Z\left( 1+\left| \tau +Dv^{\tau ,y_1}\right| +\left|
\tau +Dv^{\tau ,y_2}\right| \right) ^p\,dz\right) ^{\tfrac{\alpha (\beta -p)%
}{(\beta -1)(p-1)}}  \nonumber \\
&&\left( \int_Z\left| a(y_2,z,\tau +Dv^{\tau ,y_2})-a(y_1,z,\tau +Dv^{\tau
,y_2})\right| ^q\,dz\right) ^{\tfrac \alpha {\beta -1}}  \nonumber \\
&\leq &c\left( 1+\left| \tau \right| ^p\right) ^{\tfrac{\alpha (\beta -p)}{%
(\beta -1)(p-1)}}\widetilde{\omega }(\left| y_1-y_2\right| )\left( 1+\left|
\tau \right| ^p\right) ^{\tfrac \alpha {\beta -1}}  \nonumber \\
&\leq &\widetilde{\omega }(\left| y_1-y_2\right| )\left( 1+\left| \tau
\right| ^p\right) ^{\tfrac \alpha {p-1}}.  \label{bcontc}
\end{eqnarray}
The result follows by taking (\ref{bcont1}), (\ref{bconta}) and (\ref{bcontc}%
) into account.

(ii), (iii) and (iv): The proofs follows by similar arguments as in e.g. 
\cite{bystrom00}.%
%TCIMACRO{\TeXButton{hfill}{\hfill}}
%BeginExpansion
\hfill%
%EndExpansion
{}$\square $

\begin{remark}
\label{rema}By similar arguments it follows that (ii), (iii), and (iv) hold
up to boundaries, for the homogenized operator $b^k$ in Step 1.
\end{remark}

\label{lastpage}

\end{document}